\documentclass[leqno,12pt]{amsart}

\usepackage{amsmath,amstext,amssymb,amsopn,amsthm,mathrsfs}
\usepackage{verbatim}
\usepackage{enumerate}

\newcommand*\mr{\mathbb{R}}
\newcommand*\mn{\mathbb{N}}

\newcommand*\ven{\vert n\vert}
\newcommand*\al{\alpha}
\newcommand*\lfun{\varphi^{\alpha}_n}
\newcommand*\lfuni{\varphi^{\alpha_i}_{n_i}}
\newcommand*\lfunk{\varphi^{\alpha}_{k}}
\newcommand*\Rop{R_r^{\alpha}}
\newcommand*\Ropal{R_r^{-1/2}}
\newcommand*\Ropi{R_r^{\alpha_i}}
\newcommand*\nti{\tilde{k}}
\newcommand*\ve{\varepsilon}

\title[Hardy's inequality]
{Hardy's inequality for Laguerre expansions of Hermite type}

\author[P{.} Plewa]{Pawe\l{} Plewa}
\address{Pawe\l{} Plewa \newline
			Faculty of Pure and Applied Mathematics, 
      Wroc\l{}aw University of Science and Technology       \newline
      Wyb{.} Wyspia\'nskiego 27,
      50--370 Wroc\l{}aw, Poland      
      }
\email{pawel.plewa@pwr.edu.pl}

\allowdisplaybreaks

\theoremstyle{plain}
\newtheorem{thm}{Theorem}[section]

\newtheorem{lm}[thm]{Lemma}
\newtheorem{prop}[thm]{Proposition}
\newtheorem{cor}[thm]{Corollary}

\theoremstyle{definition}

\theoremstyle{remark}
\newtheorem*{rem*}{Remark}

\theoremstyle{plain}

\begin{document}
\begin{abstract}
Hardy's inequality for Laguerre expansions of Hermite type with the index $\al\in(\{-1/2\}\cup[1/2,\infty))^d$ is proved in the multi-dimensional setting with the exponent $3d/4$. We also obtain the sharp analogue of Hardy's inequality with $L^1$ norm replacing $H^1$ norm at the expense of increasing the exponent by an arbitrarily small value.
\end{abstract}

\maketitle
\footnotetext{
Research partially supported by funds of Faculty of Pure and Applied Mathematics, Wrocław University of Science and Technology, \#0401/0121/17.\\
The paper is a part of author's doctoral thesis written under the supervision of Professor Krzysztof Stempak.}

\section{Introduction}
The well known Hardy inequality states that
\begin{equation*}
\sum_{k\in\mathbb{Z}}\frac{\vert \hat{f}(k)\vert}{\vert k\vert+1}\lesssim \Vert f\Vert_{\mathrm{ Re} H^1}, 
\end{equation*}
where $\hat{f}(k)$ is $k$-th Fourier coefficient of $f$. Here $\mathrm{ Re} H^1$ is the real Hardy space composed of the boundary values of the real parts of functions in the Hardy space $H^1(\mathbb{D})$, where $\mathbb{D}$ is the unit disk in the plane.

Kanjin \cite{Kanjin1} established an analogue of Hardy's inequality in the context of Hermite functions $\{h_k\}_{k\in\mathbb{N}}$ and the standard Laguerre functions $\{\mathcal{L}_k^{\al}\}_{k\in\mathbb{N}}$, $\al\geq 0$, namely
\begin{equation*}
\sum_{k\in\mathbb{N}}\frac{\vert \langle f,h_k\rangle_{L^2(\mr)}\vert}{ (k+1)^{29/36}}\lesssim \Vert f\Vert_{H^1(\mr)},\qquad \sum_{k\in\mathbb{N}}\frac{\vert \langle f,\mathcal{L}_k^{\al}\rangle_{L^2(\mr_+)}\vert}{k+1}\lesssim \Vert f\Vert_{H^1(\mr_+)},
\end{equation*}
where $H^1(\mr)$ and $H^1(\mr_+)$ denote the real Hardy spaces on $\mr$ and $\mr_+$, respectively.

Hardy's inequality in the context of Hermite functions was further intensively studied by many authors. Radha \cite{Radha} proved a similar inequality in an arbitrary dimension. In \cite{RadhaThangavelu} an improved version of Hardy's inequality was introduced in the multi-dimensional case, $d\geq 2$, by Radha and Thangavelu. The exponent in the denominator was $3d/4$. This led to the hypothesis that in the one-dimensional case the exponent should be equal to $3/4$. It was indeed proved in \cite{LiYuShi} by Z. Li, Y. Yu and Y. Shi. A generalization of Kanjin's results, with the spaces $H^p(\mr)$ and $H^p(\mr_+)$, $p\in(0,1]$, instead of $H^1(\mr)$ and $H^1(\mr_+)$, was also considered in the context of Hermite functions (see \cite{BalasRadha, RadhaThangavelu}) and in the context of Laguerre functions (see \cite{RadhaThangavelu, Satake}).

In this paper we study multi-dimensional Hardy's inequality in the context of Laguerre functions of Hermite type $\{\lfun\}_{n\in\mathbb{N}^d}$. In view of the uniform boundedness of the derivatives of functions $\lfun$ and \cite[Lemma]{Kanjin1} we have the one-dimensional inequality
$$\sum_{k=0}^{\infty}\frac{\vert\langle f, \lfunk\rangle\vert}{(k+1)^{29/36}}\lesssim \Vert f\Vert_{H^1(\mr_+)}.$$
Our aim is to obtain the analogue of this inequality with the power $3d/4$, which does not depend on $\al$, and in dimension $d\geq1$.

The proof of one of the main results, Theorem \ref{main_thm}, is based on the atomic decomposition of functions from $H^1(\mr_+^d)$ and relies on a uniform estimate for atoms and an additional argument of the "weak" continuity of certain operators. Without this argument, which was often omitted in papers concerning this topic, the proof would have a gap. We remark that the uniform estimate for atoms does not imply continuity of operators that appear in analysis that involves the atomic decomposition of $H^1(\mr^d)$ (see \cite{Bownik}).

The range of the Laguerre type multi-index $\al$ that is considered in Theorem \ref{main_thm}, is the set $(\{-1/2\}\cup[1/2,\infty))^d$. This kind of restraint appeared before (see for example \cite{NowakStempak}). Note that the one-dimensional Laguerre functions of Hermite type with the Laguerre type multi-index equal to $-1/2$ or $1/2$ are, up to a multiplicative constant, the Hermite functions of even or odd degree, respectively. Therefore, it was fair to assume that this values of $\al$ should be included. Technically, the restraint emerges from the range of $\al$'s for which the derivatives of the Laguerre functions of Hermite type are uniformly bounded. It may also be related to the fact that the associated heat semi-group is a semi-group of $L^p$ contractions precisely for this set of $\al$'s (see \cite{NowakStempak2}).

In \cite{Kanjin2} Kanjin proved that if the exponent in one-dimensional Hardy's inequality in the context of Hermite functions is strictly greater than $3/4$, then one can replace $H^1(\mr)$ norm by $L^1(\mr)$ norm and, moreover, the exponent $3/4$ is sharp. In Theorem \ref{L1_thm} we shall prove that this is also the case in the context of Laguerre functions of Hermite type and extend this result to an arbitrary dimension.

We shall frequently use two basic estimates: for $a,\,A>0$ we have $\sup_{x>0} x^a e^{-Ax}< \infty$ and $(n_1+\ldots+n_d+1)^d\geq (n_1+1)\cdot\ldots\cdot(n_d+1)$, where $n_i\in\mathbb{N}$, $i=1,\ldots,d$. 

\subsection*{Notation}
Throughout this paper we write $n=(n_1,\ldots,n_d)\in\mathbb{N}^d$ for a multi-index and $\vert n\vert=n_1+\ldots+n_d$ for its length, where $\mathbb{N}=\{0,1,.\ldots\}$ and $d\geq 1$. The Laguerre type multi-index $\al=(\al_1,\ldots,\al_d)$, unless stated otherwise, is considered in the full range, i.e. $\al\in(-1,\infty)^d$. We shall also use the notation $\mr^d_+=(0,\infty)^d$ and $\mathbb{N}_+=\{1,2,\ldots\}$. For functions $f,\,g\in L^2(\mr_+^d,dx)$ we denote the inner product by $\langle f,g\rangle$. Sometimes we shall use this notation for functions that are not in $L^2(\mr_+^d,dx)$ but the underlying integral makes sense. We shall use the symbol $\lesssim$ denoting an inequality with a constant that does not depend on relevant parameters. Also, the symbol $\simeq$ means that $\lesssim$ and $\gtrsim$ hold simultaneously. Moreover, we will denote asymptotic equality by $\approx$.

\subsection*{Acknowledgement}
The author would like to thank Professor Krzysztof Stempak for insightful comments and his continuous help during the preparation of this paper.
 
\section{Preliminaries}
The {\it Laguerre functions of Hermite type} of order $\al$ on $\mr^d_+$ are the functions
$$\lfun(x)=\prod_{i=1}^d \lfuni (x_i),\qquad x=(x_1,\ldots,x_d)\in\mr^d_+, $$
where $\lfuni(x_i)$ is the one-dimensional Laguerre function of Hermite type defined by 
$$\lfuni(x_i)=\Big(\frac{2\Gamma(n_i+1)}{\Gamma(n_i+\al_i+1)}  \Big)^{1/2}L_{n_i}^{\al_i}(x_i^2)x_i^{\al_i+1/2}e^{-x_i^2/2},\qquad x_i>0.$$
The functions $\{\lfun\colon n\in\mathbb{N}^d\}$ form an orthonormal basis in $L^2(\mr_+^d,dx)$.

The one-dimensional {\it standard Laguerre functions} $\{\mathcal{L}_k^{\alpha}\}_{k\in\mathbb{N}}$ of order $\al$ are
$$\mathcal{L}_{k}^{\al}(u)=\Big(\frac{\Gamma(k+1)}{\Gamma(k+\al+1)}  \Big)^{1/2}L_{k}^{\al}(u)u^{\al/2}e^{-u/2},\qquad u>0.$$
Note that
\begin{equation*}
\lfunk(u)=(2u)^{1/2}\mathcal{L}_{k}^{\al}(u^2).
\end{equation*}
We shall use the pointwise asymptotic estimates (see \cite[p.~435]{Muckenhoupt} and \cite[p.~699]{AskeyWainger})
\begin{equation*}
\vert \mathcal{L}_{k}^{\al}(u)\vert\lesssim\left\{ \begin{array}{ll}
(u\nu)^{\al/2},  & 0<u\leq 1/\nu,\\
(u\nu)^{-1/4}, & 1/\nu<u\leq\nu/2,\\
(\nu(\nu^{1/3}+\vert u- \nu\vert))^{-1/4}, & \nu/2<u\leq 3\nu/2,\\
\exp(-\gamma u), & 3\nu/2<u<\infty,
\end{array}\right.
\end{equation*}
where $\nu=\nu(\al,k)=\max(4k+2\al+2,2)$ and with $\gamma>0$ depending only on $\alpha$. Hence,
\begin{equation}\label{pointwise_estimates_lfun}
\vert\lfunk(u)\vert\lesssim\left\{ \begin{array}{ll}
u^{\al+1/2}\nu^{\al/2},  & 0<u\leq 1/\sqrt{\nu},\\
\nu^{-1/4}, & 1/\sqrt{\nu}<u\leq\sqrt{\nu/2},\\
u^{1/2}(\nu(\nu^{1/3}+\vert u^2- \nu\vert))^{-1/4}, & \sqrt{\nu/2}<u\leq \sqrt{3\nu/2},\\
u^{1/2}\exp(-\gamma u^2), & \sqrt{3\nu/2}<u<\infty.
\end{array}\right.
\end{equation}

There is the known formula for the derivatives of functions $\lfunk$,
\begin{equation}\label{derivative_formula}
\frac{d}{d u}\lfunk(u)=-2\sqrt{k}\varphi_{k-1}^{\al+1}(u)+\left(\frac{2\al+1}{2u}-u\right)\lfunk(u),
\end{equation}
where $\varphi_{-1}^{\al}\equiv 0$.

From \eqref{pointwise_estimates_lfun} it follows that for $\al\geq -1/2$ we have
\begin{equation}\label{Linfinity_boundedness}
\Vert \lfunk\Vert_{L^{\infty}(\mr_+)}\lesssim (k+1)^{-1/12},
\end{equation}
and also by \eqref{derivative_formula}  for $\al\in\{-1,2\}\cup[1/2,\infty)$,
\begin{equation}\label{Linfinity_bound_lfun}
\left\Vert \frac{d}{d \boldsymbol{\cdot}}\lfunk(\boldsymbol{\cdot})\right\Vert_{L^{\infty}(\mr_+)}\lesssim (k+1)^{5/12}.
\end{equation}

We introduce the family of operators $\{\Rop\}_{r\in(0,1)}$, defined spectrally for $f\in L^2(\mr^d_+)$, by
\begin{equation*}
\Rop f=\sum_{n\in\mathbb{N}^d}^{\infty} r^{\vert n\vert} \langle f,\lfun\rangle\lfun.
\end{equation*}
It is easily seen by means of Parseval's identity that for every $r\in(0,1)$, the operator $\Rop$ is a contraction on $L^2(\mr_+^d)$.

The kernel associated with $\Rop$ is defined by
\begin{equation*}
\Rop(x,y)=\sum_{n\in\mathbb{N}^d} r^{\vert n\vert} \lfun(x)\lfun(y),\qquad x,y\in\mr^d_+.
\end{equation*}
Note that
$$\Rop(x,y)=\prod_{i=1}^d \Ropi(x_i,y_i) $$
and there is also the explicit formula (compare \cite[p.~102]{Szego})
\begin{equation}\label{def_R_ker_explic}
\Ropi(x_i,y_i)=\frac{2 (x_i y_i)^{1/2}}{(1-r)r^{\al_i/2}}\exp\left(-\frac{1}{2}\frac{1+r}{1-r}(x_i^2+y_i^2)\right)I_{\al_i }\left(\frac{2r^{1/2}}{1-r}x_i y_i\right),
\end{equation}
where $I_{\alpha_i}$ denotes the modified Bessel function of the first kind, which is smooth and positive on $(0,\infty)$. Notice that with $r=e^{-4t}$, $t>0$, $r^{(\vert\al\vert+d)/2}\Rop(x,y)$ is just the kernel $G_t^{\al}(x,y)$, see \cite[(2.3)]{NowakStempak}, for a differential operator $L^{\al}$ associated with $\{\lfun\}$-expansions.

Let $H^1(\mr^d)$ be the real Hardy space on $\mr^d$ (see, for example, \cite[III]{Stein}). A measurable function $a(x)$ supported in a Euclidean ball $B$ is an $H^1(\mr^d)$ atom if $\Vert a\Vert_{L^{\infty}(\mr^d)}\leq \vert B\vert^{-1}$, where $\vert B\vert$ denotes the Lebesgue measure of $B$, and $\int_B a(x)\,dx=0$. Every function $f\in H^1(\mr^d)$ has an atomic decomposition, namely there exist a sequence of complex coefficients $\{\lambda_i\}_{i=0}^{\infty}$ and a sequence of $H^1(\mr^d)$ atoms $\{a_i\}_{i=0}^{\infty}$ such that
\begin{equation}\label{atomic_decomposition}
f=\sum_{i=0}^{\infty}\lambda_i a_i, \qquad \sum_{i=0}^{\infty} \vert \lambda_i\vert\lesssim \Vert f\Vert_{H^1(\mr^d)},
\end{equation}
where the convergence of the first series is in $H^1(\mr^d)$.

The Hardy space on $\mr_+^d$ is defined by
$$H^1(\mr_+^d)=\{f\in L^1(\mr^d_+)\colon \exists \tilde{f}\in H^1(\mr^d),\ \mathrm{supp}(\tilde{f})\subset [0,\infty)^d \text{ and } \tilde{f}\big\vert_{\mr_+^d}=f\}, $$
with the norm $\Vert f\Vert_{H^1(\mr_+^d)}=\Vert \tilde{f}\Vert_{H^1(\mr^d)}$.
The properties of $H^1(\mr_+^d)$ given below follow from \cite[Lemma 7.40]{GarciaFrancia} stated in the one-dimensional case therein, however easily generalizable to the case of an arbitrary dimension. Every $f\in H^1(\mr_+^d)$ has an atomic decomposition as in \eqref{atomic_decomposition} with supports of $a_i$ in $[0,\infty)^d$; we shall call them $H^1(\mr_+^d)$ atoms. Note that a ball in $\mr_+^d$ is a ball in $\mr^d$ restricted to $\mr_+^d$. We may assume that every ball associated with an $H^1(\mr_+^d)$ atom has its center in $\mr_+^d$.

For $f\in L^1(\mr_+^d)$ we define the multi-even extension $f_e$ of $f$ by
$$f_e(x_1,\ldots,x_d)= f(\vert x_1\vert,\ldots,\vert x_d\vert),\qquad x=(x_1,\ldots,x_d)\in\mr^d.$$
We remark (again see \cite[Lemma 7.40]{GarciaFrancia}) that $f_e\in H^1(\mr^d)$ if and only if $f\in H^1(\mr_+^d)$, and $\Vert f\Vert_{H^1(\mr_+^d)}\simeq \Vert f_e\Vert_{H^1(\mr^d)}$, thus we have
\begin{equation}\label{L1&H1_inequality}
\Vert f\Vert_{L^1(\mr_+^d)}\lesssim \Vert f\Vert_{H^1(\mr_+^d)}. 
\end{equation}

\section{One-dimensional kernel estimates}\label{section_one_dimension}
We shall estimate the kernels $\Ropi(x_i,y_i)$. For the sake of convenience we will write $x,y,\al$ instead of $x_i,y_i,\al_i$.

There are known the asymptotic estimates (see \cite[p.~136]{Lebedev}) 
\begin{align*}
&I_{\al}(u) \lesssim u^{\al}, \quad 0<u<1,\\
&I_{\al}(u) \lesssim u^{-1/2}e^{u}, \quad u\geq 1.
\end{align*}
Hence, 
\begin{equation}\label{R_ker_estim}
\Rop(x,y)\lesssim\left\{\begin{array}{ll}
(1-r)^{-\al-1}(xy)^{\al+1/2}\exp\left(-\frac{1}{2}\frac{1+r}{1-r}(x^2+y^2) \right), & y\leq \frac{1-r}{2\sqrt{r}x},\\
(1-r)^{-1/2}r^{-\al/2-1/4}\exp\left(-\frac{1}{2}\frac{1+r}{1-r}(y-x)^2-\frac{xy(1-r)}{(1+\sqrt{r})^2} \right), & y\geq \frac{1-r}{2\sqrt{r}x}.
\end{array}\right.
\end{equation}

{\lm\label{lemma_R} For $\al\geq -1/2$ there is
$$\sup_{x\in\mr_+}\Vert \Rop(x,\boldsymbol{\cdot})\Vert_{L^2(\mr_+)}\lesssim (1-r)^{-1/4},\qquad r\in(0,1). $$}
\begin{proof}
For $0<r\leq 1/2$ we use Parseval's identity and \eqref{Linfinity_boundedness} obtaining
$$\sup_{x\in\mr_+}\Vert \Rop(x,\boldsymbol{\cdot})\Vert_{L^2(\mr_+)}\leq \sup_{x\in\mr_+}\Big\Vert \sum_{k\in\mathbb{N}}2^{-k} \lfunk(x)\lfunk \Big\Vert_{L^2(\mr_+)}\leq  \Big(\sum_{k\in\mathbb{N}}2^{- 2k}\Vert \lfunk\Vert_{L^{\infty}(\mr_+)}^2\Big)^{1/2}\lesssim 1.$$

For $1/2<r<1$ we denote $y_0=(1-r)/(2\sqrt{r}x)$ and estimate the integrals over $(0,y_0]$ and $(y_0,\infty)$. Thus, using the substitution $u=(y\sqrt{1+r})/\sqrt{1-r}$ we obtain
\begin{align*}
&\int_0^{y_0} \Rop(x,y)^2\,dy\\
&\lesssim \Big(\frac{x^2}{1-r}\Big)^{\al+1/2}\exp\Big(-\frac{1+r}{1-r}x^2\Big) (1-r)^{-\al-3/2}\int_0^{y_0}y^{2\al+1}\exp\Big(-\frac{1+r}{1-r}y^2\Big)\,dy\\
&\lesssim (1-r)^{-\al-3/2}\int_0^{\frac{y_0\sqrt{1+r}}{\sqrt{1-r}}}(1-r)^{\al+1/2}u^{2\al+1}e^{-u^2}(1-r)^{1/2}\,du\\
&\lesssim (1-r)^{-1/2},
\end{align*}
uniformly in $x\in\mr_+$ and $r\in(0,1)$. Similarly,
\begin{align*}
\int_{y_0}^{\infty} \Rop(x,y)^2\,dy&\lesssim  (1-r)^{-1}\int_{y_0}^{\infty}\exp\Big(-\frac{1+r}{1-r}(y-x)^2\Big)\,dy\\
&\lesssim  (1-r)^{-1}\int_{-\infty}^{\infty}\exp\Big(-\frac{1+r}{1-r}y^2\Big)\,dy\\
&\lesssim (1-r)^{-1/2},
\end{align*}
uniformly in $x\in\mr_+$ and $r\in(0,1)$. Combining the above gives the claim.
\end{proof}

{\lm\label{deri_R_ker} For $\al>0$ it holds
\begin{align*}
\partial_x\Rop(x,y)&=\Rop(x,y)\left(-\frac{2\al-1}{2x}-\frac{(1+r)x}{1-r}+\frac{2\sqrt{r}y}{1-r}\frac{I_{\al-1}(\frac{2\sqrt{r}xy}{1-r})}{I_{\al}(\frac{2\sqrt{r}xy}{1-r})} \right)\\
&=\frac{2y}{1-r}R^{\alpha-1}_r(x,y)-\left(\frac{2\al-1}{2x}+\frac{(1+r)x}{1-r}\right)\Rop(x,y). 
\end{align*}
}
\begin{proof}
It suffices to use the formula 
$$\frac{d}{du}I_{\al}(u)=-\frac{\al}{u}I_{\alpha}(u)+I_{\al-1}(u) $$
that holds for $\al>0$ (see \cite[p.~110]{Lebedev}), and  differentiate.
\end{proof}

{\lm\label{Bessel_estim} For $\alpha\geq 1/2$ there is
$$\left\vert \frac{I_{\al-1}(u)}{I_{\al}(u)}-1\right\vert\leq \frac{2\al}{u},\qquad u>0. $$ }
For the proof in the case $\al>1/2$ see \cite[pp.~6-7]{Nasell}. If $\al=1/2$, then it suffices to use the explicit formulas (see \cite[p.~112]{Lebedev})
\begin{equation}\label{explicit_formulas}
I_{-1/2}(u)=\Big(\frac{2}{\pi u}\Big)^{1/2}\cosh u,\qquad I_{1/2}(u)=\Big(\frac{2}{\pi u}\Big)^{1/2}\sinh u,\qquad u>0.
\end{equation}

Lemma \ref{Bessel_estim} is of paramount importance in our estimates wherever the cancellations are needed. It has been used before in the context of Laguerre functions (see for example \cite{NowakStempak}).

Note that Lemma \ref{deri_R_ker} works for $\al>0$, but we want to include the case $\al=-1/2$ as well. Thus, using \eqref{def_R_ker_explic} and \eqref{explicit_formulas} we obtain
$$\Ropal(x,y)=\frac{2}{\sqrt{\pi}}(1-r)^{-1/2}\exp\Big( -\frac{1}{2}\frac{1+r}{1-r}(x^2+y^2)\Big) \cosh\Big(\frac{2\sqrt{r}xy}{1-r}\Big). $$
Hence,
\begin{align}\label{ker_-1/2_deri}
\big(\partial_x\Ropal (x,y)\big)^2&=\frac{4}{\pi}(1-r)^{-3}\exp\Big(-\frac{1+r}{1-r}(x^2+y^2)\Big)\nonumber\\
&\times \bigg( 2\sqrt{r}y\sinh\Big(\frac{2\sqrt{r}xy}{1-r}\Big)-(1+r)x\cosh\Big(\frac{2\sqrt{r}xy}{1-r}\Big)\bigg)^2.
\end{align}

Using basic estimates for $\cosh$ and $\sinh$ and combining \eqref{ker_-1/2_deri} with \eqref{R_ker_estim} and Lemma \ref{deri_R_ker} we obtain for $\al\in\{-1/2\}\cup[1/2,\infty)$
\begin{align}\label{final_estimate_without_cancellations}
&\big(\partial_x\Rop (x,y)\big)^2\nonumber\\
&\lesssim\left\{\begin{array}{ll}
(1-r)^{-2\al-2}(xy)^{2\al+1}(A_\al (x,y)+x^2(1-r)^{-2})\exp\left(-\frac{1+r}{1-r}(x^2+y^2) \right), & y\leq \frac{1-r}{2\sqrt{r}x},\\
(1-r)^{-3}(x^2+y^2+(1-r)^2x^{-2})\exp\left(-\frac{1+r}{1-r}(y-x)^2-\frac{2xy(1-r)}{(1+\sqrt{r})^2} \right), & y\geq \frac{1-r}{2\sqrt{r}x},
\end{array}\right.
\end{align}
where $A_\al(x,y)=x^{-2}$ for $\al\geq 1/2$ and $A_{-1/2}(x,y)=y^2(1-r)^{-2}$.

{\prop\label{Ker_der_estim} For $\al\in\{-1/2\}\cup [1/2,\infty)$ we have
$$\sup_{x\in\mr_+}\left\Vert\partial_x\Rop(x,\boldsymbol{\cdot})\right\Vert_{L^2(\mr_+)}\lesssim (1-r)^{-3/4},\qquad r\in(0,1). $$}
\begin{proof}
Fix $x\in\mr_+$. If $0<r\leq 1/2$, then we use \eqref{Linfinity_bound_lfun} and Parseval's identity obtaining
$$\left\Vert\partial_x\Rop(x,\boldsymbol{\cdot})\right\Vert_{L^2(\mr_+)}=\Big\Vert\sum_{k\in\mathbb{N}}r^k(\lfunk)'(x)\lfunk\Big\Vert_{L^2(\mr_+)}\leq\Big(\sum_{k\in\mathbb{N}}2^{-2k}\Vert(\lfunk)'\Vert_{L^{\infty}(\mr_+)}^2\Big)^{1/2}\lesssim 1.$$ 

From now on we assume that $1/2<r<1$. We use the notation $y_0=(1-r)/(2\sqrt{r}x)$ again and split the integration over two intervals: $(0,y_0]$ and $(y_0,\infty)$. In the first case, using \eqref{final_estimate_without_cancellations} and the substitution $y=(\sqrt{1-r})/(\sqrt{1+r})u$, we obtain for $\al\geq 1/2$
\begin{align*}
&\int_0^{y_0}\left(\partial_x\Rop(x,y)\right)^2dy\\
 &\lesssim (1-r)^{-2\al-2}x^{2\al-1}\left(1+(1-r)^{-2}x^4\right)\exp{\left(-\frac{1+r}{1-r}x^2\right)}\int_0^{y_0}y^{2\al+1}\exp{\left(-\frac{1+r}{1-r}y^2\right)}dy\\
&\lesssim(1-r)^{-3/2}\left(\left(\frac{x^2}{1-r}\right)^{\al-1/2}+\left(\frac{x^2}{1-r}\right)^{\al+3/2}\right)\exp\left(-\frac{1+r}{1-r}x^2\right)\int_0^{\infty}u^{2\al+1}e^{-u^2}du\\
&\lesssim (1-r)^{-3/2}.
\end{align*}
For $\al=-1/2$ the corresponding computation is similar. The above estimate, as well as the following, are uniform in $x\in\mr_+$ and $r\in(0,1)$. 

The case of integration over $(y_0,\infty)$ is more complicated. Firstly we assume that $y_0\geq x$ and applying \eqref{final_estimate_without_cancellations} and the substitution $y-x=\sqrt{(1-r)/(1+r)}t$ we compute
\begin{align*}
\int_{y_0}^{\infty}\left(\partial_x\Rop(x,y)\right)^2 dy&\lesssim (1-r)^{-3}\int_{y_0}^{\infty}y^2\exp\left(-\frac{(1+r)(y-x)^2}{1-r}\right)dy\\
&\lesssim (1-r)^{-5/2}\int_0^{\infty}\left(t^2(1-r)+x^2\right)e^{-t^2}dt\\
&\lesssim (1-r)^{-3/2}.
\end{align*}

Now, we assume that $y_0\leq x$, and integrate over the interval $[2x,\infty)$. Similarly, we obtain
\begin{align*}
\int_{2x}^{\infty}\left(\partial_x\Rop(x,y)\right)^2 dy&\lesssim (1-r)^{-3}\int_{2x}^{\infty} y^2\exp\left(-\frac{(1+r)(y-x)^2}{1-r}\right)dy\\
&\lesssim (1-r)^{-3}\int_{x}^{\infty}(y+x)^2\exp\left(-\frac{(1+r)y^2}{1-r}\right)dy\\
&\lesssim (1-r)^{-5/2}\int_{0}^{\infty}y^2(1-r)e^{-y^2}dy\\
&\lesssim (1-r)^{-3/2}.
\end{align*}

Finally, we integrate over the interval $(y_0,2x)$ with the restrictions $1/2<r<1$ and $x\geq y_0$. Here we shall use the cancellations. Firstly we present the proof for $\al\geq 1/2$. By Lemmas \ref{deri_R_ker}, \ref{Bessel_estim}, \ref{lemma_R} and estimate \eqref{R_ker_estim} we have
\begin{align*}
&\int_{y_0}^{2x}\left(\partial_x\Rop(x,y)\right)^2dy\\
&\lesssim \int_{y_0}^{2x} \Rop(x,y)^2\Bigg(x^{-2}+\left(\frac{2\sqrt{r}y}{1-r}-\frac{(1+r)x}{1-r}\right)^2+y^2(1-r)^{-2}\bigg(1-\frac{I_{\al-1}\left(\frac{2\sqrt{r}xy}{1-r} \right)}{I_{\al}\left(\frac{2\sqrt{r}xy}{1-r}\right)}\bigg)^2\Bigg)dy\\
&\lesssim x^{-2}\int_{y_0}^{2x} \Rop(x,y)^2dy+\int_{y_0}^{2x} \Rop(x,y)^2\left(\frac{2\sqrt{r}y}{1-r}-\frac{(1+r)x}{1-r}\right)^2dy\\
&\lesssim x^{-2}(1-r)^{-1/2}\\
&\qquad+(1-r)^{-3}\int_{y_0}^{2x}\exp\Big(-\frac{(1+r)(y-x)^2}{1-r}-\frac{2xy(1-r)}{(1+\sqrt{r})^2}\Big)\big(2y\sqrt{r}-(1+r)x\big)^2dy\\
&\lesssim (1-r)^{-3/2}\\
&\qquad+(1-r)^{-3}\int_{y_0-x}^x \exp\Big(-\frac{(1+r)y^2}{1-r}-\frac{2x(y+x)(1-r)}{(1+\sqrt{r})^2}\Big)\Big(2y\sqrt{r}-\frac{x(1-r)^2}{(1+\sqrt{r})^2}\Big)^2dy\\
&\lesssim (1-r)^{-3/2}+(1-r)^{-3}\int_{y_0-x}^x y^2\exp\Big(-\frac{y^2}{1-r}\Big)dy+(1-r)x^3\exp\Big(-\frac{x^2(1-r)}{(1+\sqrt{r})^2}\Big)\\
&\lesssim (1-r)^{-3/2}+(1-r)^{-3/2}\int_{-\infty}^{\infty}y^2e^{-y^2}dy\\
&\lesssim (1-r)^{-3/2}.
\end{align*} 

Now, we consider $\al=-1/2$. We denote $z=(2\sqrt{r}xy)/(1-r)$. Equality \eqref{ker_-1/2_deri} and the estimate $\vert (1-\coth u )\sinh u \vert\leq 1$, $u>0$, yield
\begin{align*}
&\big(\partial_x\Ropal (x,y)\big)^2\\
&\qquad\lesssim(1-r)^{-3}\exp\Big(-\frac{1+r}{1-r}(x^2+y^2)\Big) \big( 2\sqrt{r}y-(1+r)x\coth z \big)^2\sinh^2 z\\
&\qquad\lesssim (1-r)^{-3}\exp\Big(-\frac{1+r}{1-r}(x^2+y^2)\Big) \Big(x^2+\big( 2\sqrt{r}y-(1+r)x\big)^2\sinh^2 z\Big).
\end{align*}
Note that
\begin{align*}
&(1-r)^{-3}x^2\int_{y_0}^{2x}\exp\Big(-\frac{1+r}{1-r}(x^2+y^2)\Big)\,dy\\
&\lesssim (1-r)^{-2}\frac{x^2}{1-r}\exp\Big(-\frac{1+r}{1-r}x^2\Big)\int_{y_0}^{2x}\exp\Big(-\frac{1+r}{1-r}y^2\Big)\,dy\\
&\lesssim (1-r)^{-3/2}.
\end{align*}
Moreover, using the estimate for the hyperbolic sine we obtain
\begin{align*}
&(1-r)^{-3}\int_{y_0}^{2x}\exp\Big(-\frac{1+r}{1-r}(x^2+y^2)\Big)\sinh^2\Big(\frac{2\sqrt{r}xy}{1-r}\Big)\big( 2\sqrt{r}y-(1+r)x\big)^2\,dy\\
&\lesssim (1-r)^{-3}\int_{y_0}^{2x}\exp\Big(-\frac{1+r}{1-r}(y-x)^2-\frac{2xy(1-r)}{(1+\sqrt{r})^2}\Big)\big( 2\sqrt{r}y-(1+r)x\big)^2\,dy,
\end{align*}
but this is the same quantity as in the corresponding estimate in the case $\al\geq 1/2$.
\end{proof}

Now we can state the multi-dimensional corollary.

{\cor\label{cor_of_prop} For $\al\in (\{-1/2\}\cup[1/2,\infty))^d$ and $j\in\{1,\ldots,d\},$ we have
$$\sup_{x\in\mr_+^d}\left\Vert\partial_{x_j}\Rop(x,\boldsymbol{\cdot})\right\Vert_{L^2(\mr_+^d)}\lesssim (1-r)^{-(d+2)/4}, \qquad\ r\in(0,1). $$ }
\begin{proof}
For simplicity we can assume that $j=1$. Thus,
$$\partial_{x_1}\Rop(x,y)=\frac{\partial}{\partial x_1}R^{\alpha_1}_r(x_1,y_1)\prod_{i=2}^d \Ropi(x_i,y_i). $$
Hence, Lemma \ref{lemma_R} and Proposition \ref{Ker_der_estim} imply
$$\left\Vert\partial_{x_1}\Rop(x,\boldsymbol{\cdot})\right\Vert_{L^2(\mr_+^d)}=\left\Vert \partial_{x_1}R^{\alpha_1}_r(x_1,\boldsymbol{\cdot})\right\Vert_{L^2(\mr_+)}\prod_{i=2}^d \left\Vert \Ropi(x_i,\boldsymbol{\cdot})\right\Vert_{L^2(\mr_+)}\lesssim (1-r)^{-(d+2)/4},$$
uniformly in $x=(x_1,\ldots,x_d)\in\mr_+^d$ and $r\in(0,1)$.
\end{proof}
\section{Main results}

{\prop\label{prop_atoms} For $\al\in (\{-1/2\}\cup[1/2,\infty))^d$ there is
$$\int_0^1 \Vert \Rop a\Vert_{L^2(\mr_+^d)}(1-r)^{(d-4)/4}dr\lesssim 1,  $$
uniformly in $H^1(\mr_+^d)$ atoms $a$.}
\begin{proof}
Let us fix an $H^1(\mr_+^d)$ atom $a$ supported in a ball $B$. Let $x'=(x_1',\ldots,x_d')\in\mr_+^d$ be the center of $B$. Note that since $\Rop$ are contractions on $L^2(\mr_+^d)$ we have for every $0<r<1$
$$\Vert \Rop a\Vert_{L^2(\mr_+^d)}\leq\Vert a\Vert_{L^2(\mr_+^d)}\leq \vert B\vert^{-1/2}.$$
This finishes the proof in case $\vert B\vert\geq 1$. From now on, let us assume $\vert B\vert<1$. Minkowski's integral inequality and Corollary \ref{cor_of_prop} imply 
\begin{align*}
\Vert \Rop a&\Vert_{L^2(\mr_+^d)}\\
&=\Big( \int_{\mr_+^d}\Big\vert\int_{B}\big(\Rop(x_1,x_2,\ldots,x_d,y)-\Rop(x_1',\ldots,x_d',y)\big)a(x)dx\Big\vert^2  dy\Big)^{1/2}\\
&=\Big( \int_{\mr_+^d}\Big\vert\int_{B}\Big(\sum_{i=1}^d \int_{x_i'}^{x_i} \partial_{x_i}\Rop(x'_1,\ldots,x'_{i-1},s,x_{i+1},\ldots,x_d,y)\,ds\Big)a(x)\,dx\Big\vert^2  dy\Big)^{1/2}\\
&\lesssim \int_{B} \vert a(x)\vert \vert B\vert^{1/d}\sum\limits_{i=1}^d\sup_{\xi\in\mr_+^d}\big\Vert\partial_{x_i}\Rop(\xi,\boldsymbol{\cdot})\big\Vert_{L^2(\mr_+^d)}\, dx\\ 
&\lesssim \vert B\vert^{1/d} (1-r)^{-(d+2)/4}.
\end{align*}
Thus, using the above estimates we obtain
\begin{align*}
\int_0^1 \Vert &\Rop a\Vert_{L^2(\mr^d_+)}(1-r)^{(d-4)/4}dr\\
&\lesssim \int_0^{1-\vert B\vert^{2/d}} \vert B\vert^{1/d} (1-r)^{-3/2}dr+\int_{1-\vert B \vert^{2/d}}^1 \vert B\vert^{-1/2}(1-r)^{(d-4)/4}dr,
\end{align*}
and this quantity is bounded by a constant that does not depend on $\vert B\vert$.
\end{proof}
Now we can state the main theorem.

{\thm\label{main_thm} Let $\al\in (\{-1/2\}\cup[1/2,\infty))^d$. Then
$$\sum_{n\in\mathbb{N}^d}\frac{\vert\langle f,\lfun\rangle\vert}{(\vert n\vert+1)^{3d/4}}\lesssim  \Vert f\Vert_{H^1(\mr_+^d)},$$
uniformly in $f\in H^1(\mr_+^d)$.}

\begin{proof}
Firstly we prove that
$$\sum_{n\in\mathbb{N}^d}\frac{\vert\langle a,\lfun\rangle\vert}{(\vert n\vert+1)^{3d/4}}\lesssim 1,$$
uniformly in $H^1(\mr_+^d)$ atoms $a$.

We shall employ the same argument that is used in \cite{LiYuShi}. For the Beta function there is the known asymptotic $B(k,m)\approx\Gamma(m)k^{-m}$ for large $k$ and fixed $m$. Let $a$ be an $H^1(\mr_+^d)$ atom. Using H\"{o}lder's inequality and Proposition \ref{prop_atoms} we obtain
\begin{align*}
\sum_{n\in\mathbb{N}^d}\frac{\vert\langle  a, \lfun\rangle\vert}{(\vert n\vert+1)^{3d/4}}&\lesssim\sum_{n\in\mathbb{N}^d}\int_0^1 r^{2\vert n\vert}(1-r)^{(3d-4)/4}\vert\langle a, \lfun\rangle\vert dr\\
&\leq\int_0^1 (1-r)^{(3d-4)/4}\left(\sum_{n\in\mathbb{N}^d}r^{2\vert n\vert}\right)^{1/2}\left(\sum_{n\in\mathbb{N}^d}r^{2\vert n\vert}\vert\langle a,\lfun\rangle\vert^2\right)^{1/2}dr\\
&\lesssim\int_0^1 (1-r)^{(3d-4)/4}(1-r)^{-d/2}\Vert \Rop a\Vert_{L^2(\mr_+^d)}dr\\
&\lesssim 1.
\end{align*}

Now, we define $T(f)=\{\langle f,\lfun\rangle\}_{n\in\mathbb{N}^d}$ for $f\in H^1(\mr_+^d)$. Our aim is to prove that $T\colon H^1(\mr_+^d)\rightarrow \ell^1((\vert n\vert+1)^{-3d/4})$, is bounded. Note that \eqref{Linfinity_boundedness} and \eqref{L1&H1_inequality} yield
$$\vert\langle f,\lfun\rangle\vert\leq \Vert \lfun\Vert_{L^{\infty}(\mr_+^d)} \Vert f\Vert_{L^1(\mr_+^d)}\lesssim\prod_{i=1}^d (n_i+1)^{-1/12}\Vert f\Vert_{H^1(\mr_+^d)}\leq (\vert n\vert +1)^{-1/12}\Vert f\Vert_{H^1(\mr_+^d)}.$$
Thus, $T\colon H^1(\mr_+^d)\rightarrow \ell^1((\vert n\vert+1)^{-d})$ is bounded. Note also that
\begin{equation}\label{ell_norms_ineq}
 \Vert \boldsymbol{\cdot}\Vert_{\ell^1((\vert n\vert+1)^{-d})}\leq \Vert \boldsymbol{\cdot}\Vert_{\ell^1((\vert n\vert+1)^{-3d/4})}.
\end{equation}

Let us take $f\in H^1(\mr_+^d)$ and $f=\sum_{i=0}^{\infty}\lambda_i a_i$ be an atomic decomposition of $f$. Denote $f_m=\sum_{i=0}^{m}\lambda_i a_i$ and note that $T(f_m)$ is a Cauchy sequence in $\ell^1((\vert n\vert+1)^{-3d/4})$. Indeed, we have for $l<m$,
$$\Vert T(f_m)-T(f_l) \Vert_{\ell^1((\vert n\vert+1)^{-3d/4})}\leq \sum_{i=l+1}^{m}\vert\lambda_i\vert\Vert T(a_i)\Vert_{\ell^1((\vert n\vert+1)^{-3d/4})}\lesssim \sum_{i=l+1}^m \vert \lambda_i\vert. $$
Hence, $T(f_m)$ converges to a sequence $g$ in $\ell^1((\vert n\vert+1)^{-3d/4})$ and, by \eqref{ell_norms_ineq}, also in $\ell^1((\vert n\vert+1)^{-d})$. Since $T\colon H^1(\mr_+^d)\rightarrow \ell^1((\vert n\vert+1)^{-d})$ is bounded we have $T(f_m)\rightarrow T(f)$ in $\ell^1((\vert n\vert+1)^{-d})$, therefore $g=T(f)$. To obtain the boundedness of $T\colon H^1(\mr_+^d)\rightarrow \ell^1((\vert n\vert+1)^{-3d/4})$ we fix $\ve>0$ and take $m$ such that $\Vert T(f-f_m)\Vert_{\ell^1((\ven+1)^{-3d/4})}<\ve$ and calculate
\begin{align*}
\Vert T(f)\Vert_{\ell^1((\ven+1)^{-3d/4})}&\leq \Vert T(f-f_m)\Vert_{\ell^1((\ven+1)^{-3d/4})}+\Vert T(f_m)\Vert_{\ell^1((\ven+1)^{-3d/4})}\\
&\leq \ve +\sum_{i=0}^m \vert \lambda_i\vert \Vert T(a_i)\Vert_{\ell^1((\ven+1)^{-3d/4})}\\
&\lesssim \ve+\Vert f\Vert_{H^1(\mr_+^d)}.
\end{align*}
This finishes the proof.
\end{proof}

\section{$L^1$ result}
In this section we shall prove that the inequality in Theorem \ref{main_thm} holds also with $L^1(\mr_+^d)$ norm replacing $H^1(\mr_+^d)$ norm provided that the exponent in the denominator is strictly greater than $3d/4$. Our reasoning is similar to Kanjin's in \cite{Kanjin2}. The main tool in the proof of this fact is the asymptotic estimate for functions $\lfun$.
{\thm\label{L1_thm} Let $\ve>0$ and $\al\in [-1/2,\infty)^d$. Then 
\begin{equation}\label{L1_thm_ineq}
\sum_{n\in\mathbb{N}^d}\frac{\vert\langle f,\lfun\rangle\vert}{(\vert n\vert+1)^{\frac{3d}{4}+\ve}}\lesssim \Vert f\Vert_{L^1(\mr_+^d)},
\end{equation}
uniformly in $f\in L^1(\mr_+^d)$.
The result is sharp in the sense that there is $f\in L^1(\mr_+^d)$ such that
\begin{equation}\label{L1_counterexample}
\sum_{n\in\mathbb{N}^d}\frac{\vert\langle f,\lfun\rangle\vert}{(\vert n\vert+1)^{3d/4}}=\infty.
\end{equation}
} 
\begin{proof}
Given $\ve>0$ and $\al\in[-1/2,\infty)^d$, for the proof of \eqref{L1_thm_ineq} it suffices to verify that
$$\sum_{n\in\mathbb{N}^d}\frac{\vert\lfun(x)\vert}{(\vert n\vert+1)^{\frac{3d}{4}+\ve}}\lesssim 1,\qquad x\in\mr_+^d.$$
We shall prove this estimate in the one-dimensional case. This is indeed sufficient, since
$$\sum_{n\in\mathbb{N}^d}\frac{\vert\lfun(x)\vert}{(\vert n\vert+1)^{\frac{3d}{4}+\ve}}\leq \prod_{i=1}^d \sum_{n_i=0}^{\infty}\frac{\vert\lfuni(x_i)\vert}{( n_i+1)^{\frac{3}{4}+\ve/d}}. $$
But $\varphi^{\al}_0(u)\lesssim 1$ uniformly in $u\in\mr_+$, so given $\al\in[-1/2,\infty)$ we are reduced to proving
\begin{equation}\label{L1_new_claim}
\sum_{k=1}^{\infty}\frac{\vert \lfunk(u)\vert}{k^{3/4+\ve}}\lesssim 1,\qquad u\in\mr_+.
\end{equation}

Denote $\nti=4k+2\al+2$ and for $u\in\mr_+$ define 
$$\mathcal{N}_u=\big\{k\in\mathbb{N}_+\colon \nti/2\leq u^2\leq 3\nti/2\big\}.$$
We have
\begin{equation*}
\sum_{k=1}^{\infty}\frac{\vert \lfunk(u)\vert}{k^{3/4+\ve}}=\sum_{k\notin \mathcal{N}_u}\frac{\vert \lfunk(u)\vert}{k^{3/4+\ve}}+\sum_{k\in \mathcal{N}_u}\frac{\vert \lfunk(u)\vert}{k^{3/4+\ve}}.
\end{equation*}
Note that by \eqref{pointwise_estimates_lfun}, if $k\notin \mathcal{N}_u$, then $\vert \lfunk(u)\vert\lesssim k^{-1/4}$ uniformly in $u$ and $k$, and hence the sum over the complement of $\mathcal{N}_u$ is bounded uniformly in $u\in\mr_+$. We claim that the same is true for the sum over $\mathcal{N}_u$.  

Assume $\mathcal{N}_u\neq \emptyset$ and let $k_0=k_0(u)=\min\{k\in\mathbb{N}_+\colon k\in \mathcal{N}_u\}$. Definition of $\mathcal{N}_u$ implies that $\mathcal{N}_u\subset [k_0,k^*_0]$, where $k^*_0=3k_0+1+\lceil\al\rceil$. Thus, \eqref{pointwise_estimates_lfun} yields
\begin{align*}
\sum_{k\in \mathcal{N}_u}\frac{\vert \lfunk(u)\vert}{k^{3/4+\ve}}\lesssim \sum_{k=k_0}^{k^*_0}\frac{\sqrt{u}}{k^{3/4+\ve}\,\nti^{1/4}(\nti^{1/3}+\vert u^2-\nti\vert)^{1/4}}&\lesssim \sum_{k=k_0}^{k^*_0}\frac{(3\widetilde{k_0}/2)^{1/4}}{k_0^{1+\ve}(\nti^{1/3}+\vert u^2-\nti\vert)^{1/4}}\\
&\lesssim k_0^{-3/4-\ve}\sum_{k=k_0}^{k^*_0}\frac{1}{(1+\vert u^2-\nti\vert)^{1/4}},
\end{align*}
uniformly in $u$. Note that $\vert u^2-\nti\vert$ increases in $k$ provided $\nti\geq u^2$. Since $u^2\leq 3\widetilde{k_0}/2$, we have for $\nti\geq 3\widetilde{k_0}/2$ or, equivalently, for $k\geq k_0^{**}:=\lceil 3k_0/2+\al/4+1/4\rceil$, that $\vert u^2-\nti\vert$ increases in $k$. Hence,
\begin{align*}
\sum_{k=k_0}^{k^*_0}\frac{1}{(1+\vert u^2-\nti\vert)^{1/4}}\lesssim\sum_{k=k_0^{**}}^{k^*_0}\frac{1}{(1+\vert u^2-\nti\vert)^{1/4}}&\lesssim \int_{k_0^{**}}^{k^*_0}\frac{dt}{(1+4t+2\al+2-u^2)^{1/4}}\\
&\lesssim (4k^*_0+2\al+3-u^2)^{3/4}\\
&\simeq k_0^{3/4},
\end{align*}
uniformly in $u$. This completes the proof of the claim and hence the justification of \eqref{L1_new_claim} and thus finishes the verification of \eqref{L1_thm_ineq}.

Now we pass to the proof of \eqref{L1_counterexample}. Assume a contrario that the sum in \eqref{L1_counterexample} is finite for every $f\in L^1(\mr_+^d)$. The uniform boundedness principle implies that
$$\sum_{n\in\mathbb{N}^d}\frac{\vert \langle f,\lfun\rangle\vert}{(\vert n\vert+1)^{3d/4}}\lesssim \Vert f\Vert_{L^1(\mr_+^d)}, $$
uniformly in $f\in L^1(\mr_+^d)$. Hence, by an obvious adaptation of \cite[Lemma~1]{Kanjin2} we obtain
\begin{equation*}
\sum_{n\in\mathbb{N}^d}\frac{\vert \lfun(x)\vert}{(\vert n\vert+1)^{3d/4}}\lesssim 1,\qquad x\in\mr_+^d.
\end{equation*}

But, as we shall see, it does not hold. In fact, we shall prove that for any $x\in\mr_+^d$ we have
\begin{equation}\label{claim_L1_counterexample}
\sum_{n\in\mn^d_+}\frac{\vert \lfun(x)\vert}{\ven^{3d/4}}=\infty.
\end{equation}

Notice that using the asymptotic estimate for Laguerre polynomials (see \cite[(4.22.19)]{Lebedev}) and the known asymptotic for the Gamma function, $\Gamma(k+a)/\Gamma(k+b)\approx k^{a-b}$, $k\rightarrow\infty$, where $a,\,b\geq 0$ are fixed, we obtain for $u\in\mr_+$ and $\beta\geq-1/2$
$$\varphi_k^{\beta}(u)\approx  \pi^{-1/2} k^{-1/4}\cos\Big(2\sqrt{k}u-\frac{\pi(2\beta+1)}{4}\Big),\qquad k\rightarrow \infty. $$
Hence, we reduce verifying \eqref{claim_L1_counterexample} to checking that
\begin{equation}\label{claim_cosinus}
\sum_{n\in\mn_+^d}\frac{\bigg\vert \prod_{i=1}^d \cos\Big(2\sqrt{n_i}x_i-\frac{\pi(2\al_i+1)}{4}\Big)\bigg\vert}{\ven^d}=\infty.
\end{equation}

We first prove the one-dimensional case. Fix $u\in\mr_+$ and notice that, for $d=1$, the corresponding sum in \eqref{claim_cosinus} is greater than
\begin{align*}
\sum_{k=1}^{\infty}\frac{ \cos^2\Big(2\sqrt{k}u-\frac{\pi(2\beta+1)}{4}\Big)}{k}=\sum_{k=1}^{\infty}\frac{ 1+\cos 4\sqrt{k}u\cos\frac{\pi(2\beta+1)}{2}+\sin4\sqrt{k}u\sin\frac{\pi(2\beta+1)}{2}}{2k}.
\end{align*}
Thus, \eqref{claim_cosinus} holds, since for any $t\in\mr\smallsetminus\{0\}$ each of the two series
\begin{equation*}
\sum_{k=1}^{\infty}\frac{ \Big\{\substack{\sin\\ \cos}\Big\}\, \big(t\sqrt{k}\big)}{k}
\end{equation*}
converges. Since we could not find a proof of this fact in the literature, we offer a short argument (for the cosine series and $t=1$).

Let $H(k)=\sum_{j=1}^{k}1/j$ denote the $k$-th harmonic number. Applying  summation by parts, for any $K\in\mathbb{N}_+$ we obtain
\begin{equation}\nonumber
\sum_{k=1}^{K}\frac{\cos\sqrt{k}}{k}=H(K)\cos\sqrt{K}+\int_1^K H(\lfloor u\rfloor)\frac{\sin\sqrt{u}}{2\sqrt{u}}\,du.
\end{equation}
We use the asymptotic $H(k)=\log k+\gamma+r(k)$, where $r(k)=O(1/k)$ and $\gamma$ is the Euler-Mascheroni constant, and plug it into the both summands on the right hand side of the above formula. The terms resulting from the error parts, namely $r(K)\cos\sqrt{K}$ and 
$$\int_1^K r(\lfloor u\rfloor)\frac{\sin\sqrt{u}}{2\sqrt{u}}\,du$$
are easily seen to converge with $K\rightarrow\infty$. This is also true for
$$\int_1^{K}\log\lfloor u\rfloor\frac{\sin\sqrt{u}}{\sqrt{u}}\,du-\int_1^{K}\log u\,\frac{\sin\sqrt{u}}{\sqrt{u}}\,du=\int_1^{K}\log\Big(1+\frac{\lfloor u\rfloor-u}{u}\Big)\frac{\sin\sqrt{u}}{\sqrt{u}}\,du.$$
Thus we are left with
\begin{align*}
(\log K+\gamma)\cos\sqrt{K}+\int_1^K (\log u+\gamma)\frac{\sin\sqrt{u}}{2\sqrt{u}}\,du&=\gamma\cos1+\int_1^K \frac{\cos\sqrt{u}}{u}\,du.
\end{align*}
The latter integral, after a change of variable, is also easily seen to converge with $K\rightarrow\infty$. This finishes the proof of the convergence of the investigated series.

Now we continue and prove \eqref{claim_cosinus} in the multi-dimensional setting. Given $x\in\mr_+^d$ and proceeding similarly as before we reduce justifying \eqref{claim_cosinus} to verifying that each of the $3^d-1$ iterated series
\begin{equation}\label{all_series}
\sum_{n_1=1}^{\infty}\ldots\sum_{n_d=1}^{\infty}\vert n\vert^{-d}\prod_{j\in J}\Big\{\substack{\sin\\ \cos}\Big\}\big(t_j\sqrt{n_j}\big)
\end{equation}
converges, where $J$ is any non-empty subset of $\{1,\ldots,d\}$ and $t_j\neq0$, $j\in J$.  We shall use the induction over the dimension. Suppose that every series of the form as in \eqref{all_series} converges. We will prove that also the analogous series in dimension $d+1$ converge. Fix such a series and consider the associated set $J\subset\{1,\ldots,d+1\}$. We distinguish two cases depending on whether $d+1\in J$ or not.

If $d+1\notin J$, then the investigated series is of the form
$$\sum_{n_1=1}^{\infty}\ldots\sum_{n_d=1}^{\infty}\prod_{j\in J}\Big\{\substack{\sin\\ \cos}\Big\}\big(t_j\sqrt{n_j}\big)\sum_{k=1}^{\infty}(\vert n\vert+k)^{-d-1}.$$
It now suffices to use the asymptotic
$$\sum_{k=1}^{\infty}(\vert n\vert +k)^{-d-1}=\vert n\vert^{-d}+O(\vert n\vert^{-d-1})$$
and the inductive assumption.

The case $d+1\in J$ is more involved. We simplify matters, without any loss o generality, assuming $t_j=1$. The considered series is of the form 
$$\sum_{n_1=1}^{\infty}\ldots\sum_{n_d=1}^{\infty}\Lambda(J,n)\sum_{k=1}^{\infty}\frac{\cos\sqrt{k}}{(\vert n\vert+k)^{d+1}}$$
(or with the sine in place of the cosine, but this is not an obstacle), where $\Lambda(J,n)$ is a product of the sines or the cosines taken at $\sqrt{n_j}$, $j\in J$, respectively. In fact, we shall prove the slightly stronger result that
\begin{equation}\label{abs_convergence}
\sum_{n\in\mn^d_+}\Big\vert\sum_{k=1}^{\infty}\frac{\cos\sqrt{k}}{(\vert n\vert+k)^{d+1}}\Big\vert<\infty.
\end{equation}
We remark that the cancellation provided by one trigonometric functions are sufficient in our estimates. Note that we cannot use the triangle inequality in the innermost series, because the resulting series would diverge.

To verify \eqref{abs_convergence} we check the convergence of the innermost series with a control of the decrease of its sum in $\vert n\vert$. We will use the following asymptotic estimate
$$\sum_{k=1}^{\lfloor u\rfloor}\frac{1}{\ven +k}=\log\Big(1+\frac{u}{\ven} \Big)+r_u(\vert n\vert),$$
where $r_u(\ven)=O(\ven^{-1})$ uniformly in $u\in\mr_+$. Summation by parts and the above asymptotic yield
\begin{align*}
&\sum_{k=1}^{K}\frac{1}{\ven +k}\frac{\cos\sqrt{k}}{(\ven+k)^{d}}\\
&=\left(\sum_{k=1}^{K}\frac{1}{\ven +k}\right)\frac{\cos\sqrt{K}}{(\ven+K)^{d}}-\int_1^K\left(\sum_{k=1}^{\lfloor u\rfloor}\frac{1}{\ven +k}\right)\bigg(\frac{\cos\sqrt{u}}{(\ven +u)^{d}}\bigg)'\,du\\
&=\log\Big(1+\frac{K}{\ven}\Big)\frac{\cos\sqrt{K}}{(\ven+K)^{d}}-\int_1^K\log\Big(1+\frac{u}{\ven}\Big)\bigg(\frac{\cos\sqrt{u}}{(\ven+u)^{d}}\bigg)'\,du\\
&+ r_K(\ven) \frac{\cos\sqrt{K}}{(\ven+K)^{d}}-\int_1^K r_u(\ven)\bigg(\frac{\cos\sqrt{u}}{(\ven+u)^{d}}\bigg)'\,du.
\end{align*}
The term with the error part $r_K(\ven)$ converges to zero with $K\rightarrow\infty$, while the integral term of the error part $r_u(\ven)$ is absolutely convergent with proper decrease in $\ven$, namely
$$\int_1^K \Big\vert r_u(\ven)\bigg(\frac{\cos\sqrt{u}}{(\ven+u)^{d}}\bigg)'\Big\vert\,du\lesssim \ven^{-1}\int_1^{\infty} \frac{1}{\sqrt{u}(\ven+u)^d}\,du\lesssim \ven^{-(d+1/4)}. $$
On the other hand, for the main terms, using integration by parts twice we obtain
\begin{align*}
&\log\Big(1+\frac{K}{\ven}\Big)\frac{\cos\sqrt{K}}{(\ven+K)^{d}}-\int_1^K\log\Big(1+\frac{u}{\ven}\Big)\bigg(\frac{\cos\sqrt{u}}{(\ven+u)^{d}}\bigg)'\,du\\
&=\log\Big(1+\frac{1}{\ven}\Big)\frac{\cos1}{(\ven+1)^{d}}+\int_1^K\frac{\cos\sqrt{u}}{(\ven+u)^{d+1}}\,du\\
&=\log\Big(1+\frac{1}{\ven}\Big)\frac{\cos1}{(\ven+1)^{d}}+\frac{2\sqrt{K}\sin\sqrt{K}}{(\ven+K)^{d+1}}-\frac{2\sin1}{(\ven+1)^{d+1}}\\
&\qquad+\int_1^K \frac{\sin\sqrt{u}\,((2d+1)u-\ven)}{\sqrt{u}(\ven+u)^{d+2}}\,du.
\end{align*}
Thus, combining the above and passing to the limit with $K\rightarrow\infty$ we get
\begin{align*}
\Big\vert\sum_{k=1}^{\infty}\frac{\cos\sqrt{k}}{(\ven+k)^{d+1}}\Big\vert&\lesssim \ven^{-d-1}+\int_1^{\infty} \frac{(2d+1)u+\ven}{\sqrt{u}(\ven+u)^{d+2}}\,du+\ven^{-d-1/4}\lesssim \ven^{-d-1/4}.
\end{align*}
This finishes the justification of the convergence of the considered series and thus the verification of \eqref{claim_cosinus}. The validation of \eqref{L1_counterexample} is completed and also the proof of the whole theorem is finished.
\end{proof}

\end{document}